\documentclass[12pt,leqno]{article}
\usepackage{latexsym}
\usepackage{bbm}
\usepackage{amssymb,amsmath}
\begin{titlepage}
\title{Addendum to ``On the $p^{\lambda}$ problem''}
\author{Stephan Baier
\thanks{e-mail: sbaier@mast.queensu.ca}}
\end{titlepage}
\begin{document}
\maketitle
{\bf Abstract:} We prove that the conditions $\lambda<5/19$ and 
$L\le T^{1/2}$ in Theorems 3 and 4 of our recent paper [Bai] can be
omitted.\\ \\

In [Bai] we proved the following mean value estimate for products of
shifted and ordinary Dirichlet polynomials.\\ 

{\sc Theorem 1:} ([Bai, Theorem 4]) \begin{it} 
Suppose that  $\alpha\not=0$, $0\le \theta<1$, $T> 0$, 
$K\ge 1$, $L\ge 1$. If $\theta\not=0$, then additionally suppose that 
$L\le T^{1/2}$. 
Let  $(a_{k})$ and  $(b_{l})$ be arbitrary sequences of complex
numbers. Suppose that $\vert a_k \vert \le A$ for all $k\sim K$ and 
$\vert b_l \vert \le B$ for all $l\sim L$. Then,
  \begin{eqnarray}
  & & \int\limits_0^T \left\vert \sum\limits_{k\sim K} 
  a_{k}k^{it} \right\vert^2 \
  \left\vert \sum\limits_{l\sim L} 
  b_{l}(l+\theta)^{i\alpha t} \right\vert^2  {\rm d}t \label{15} \\ 
  \nonumber\\
  &\ll& A^2B^2(T+KL)KL\log^3(2KLT), \nonumber
  \end{eqnarray}
the implied $\ll$-constant depending only on $\alpha$. If $\theta=0$, then
$\log^3(2KLT)$ on the right side of {\rm(\ref{15})} may be replaced by 
$\log^2(2KLT)$.\end{it}\\

We then used this mean value estimate to prove the following result on the
$p^{\lambda}$ problem.\\

{\sc Theorem 2:} ([Bai, Theorem 3]) \begin{it} 
Suppose that $\varepsilon>0$, $B>0$, $\lambda\in (0,1/2]$ 
and a real $\theta$ are given. If $\theta$ is irrational, then suppose 
that $\lambda< 5/19$. Let $N\ge 3$. 
Let ${\bf A}$ be an arbitrarily given subset of
the set of positive integers. Define
  \begin{displaymath}
  F_{\theta}(\lambda):=\left\{ \begin{array}{llll} 
  \max\limits_{k\in {\bf N}}\
  \min\left\{\frac{5}{12}-\frac{(k+6)\lambda}{6(k+1)},
  \frac{5}{11}-\frac{(5k+1)\lambda}{11}\right\}& 
  \mbox{ if } \theta \mbox{ is 
  rational,} \\ \\ \frac{5}{12}-\frac{7\lambda}{12} 
  & \mbox{ otherwise.} \end{array} \right.
  \end{displaymath}
Suppose that
  $$
  N^{-F_{\theta}(\lambda)+\varepsilon\lambda}\le \delta \le 1. 
  $$
Then we have 
  $$
  \sum\limits_{\scriptsize \begin{array}{cccc} &N<n\le 2N,&\\ 
  &\left\{n^{\lambda}-\theta\right\}<\delta,&\\ 
  &\left[n^{\lambda}\right]\in {\bf A}& \end{array}} \hspace{-0.4cm}
  \Lambda(n)\ =\ \frac{\delta}{\lambda} \cdot\hspace{-0.4cm}
  \sum\limits_{\scriptsize \begin{array}{cccc} 
  &N^{\lambda}<n\le (2N)^{\lambda},&\\ 
  &n\in {\bf A}& \end{array}} \hspace{-0.4cm}
  n^{1/\lambda-1}\ + \
  O\left(\frac{\delta N}{(\log N)^B}\right). 
  $$
\end{it}
 
At the end of the last section in [Bai] we pointed out that if the condition
$L\le T^{1/2}$ in the above Theorem 1 could be omitted, then the condition 
$\lambda<5/19$ in Theorem 2 could be omitted too. 
In the following we will see 
that the condition $L\le T^{1/2}$ in Theorem 1
is actually superfluous if we allow ourselves to weaken the mean value estimate
(\ref{15}) slightly. We establish the following\\ 

{\sc Theorem 3:} \begin{it} Let 
$\theta,\xi,\alpha,\beta$
be real numbers with $0\le \theta,\xi<1$ and 
$\alpha\beta\not=0$. Suppose that $T, K, L \ge 1$, 
$\vert a_k\vert \le 1$ and $\vert b_l\vert \le 1$. Then
\begin{equation}
\int\limits_{0}^{T} \left\vert \sum\limits_{k\sim K} a_{k}
(k+\theta)^{i\alpha t} \right\vert^{2} 
\left\vert \sum\limits_{l\sim L} b_{l}
(l+\xi)^{i\beta t} \right\vert^{2} {\rm d}t
\ll \left(T+KL\right)KL
(\log T)^{15}. \label{100} 
\end{equation}
\end{it}

In accordance with the proof of [Bai, Theorem 3], from the above 
Theo\-rem 3 with
$\xi=0$ it can be deduced that Theorem 2 holds true with the condition 
$\lambda<5/19$ omitted.
    
The main idea of our proof of Theorem 3 is to relate the shifted Dirichlet
polynomials on the left-hand side of (\ref{100}) to the corresponding
Hurwitz zeta functions. For technical reasons we here define the Hurwitz
zeta function $\zeta(s,y)$ in a slightly different manner to normal usage.
For $0\le y< 1$ and Re $s>1$ we write
$$
\zeta(s,y):=\sum\limits_{n=1}^{\infty} (n+y)^{-s}.
$$
In the usual definition the series on the right-hand
side starts with $n=0$, and the
case $y=0$ is excluded, which we seek to avoid here.
 
As a function of $s$, the Hurwitz zeta function  
has a meromorphic conti\-nuation to the entire complex
plane, with a simple pole at $s=1$ (see [Ivi]).
At first, we establish the following fourth power moment estimate for the
Hurwitz zeta function on the critical line.\\

{\sc Theorem 4:} \begin{it} Suppose that $V>2\pi$ and $0\le y<1$. Then
$$  
\int\limits_{-V}^{V} \left\vert \zeta\left(\frac{1}{2}+it,y\right)
\right\vert^{4} {\rm d}t
\ll V(\log V)^{10}.
$$\end{it}

{\sc Proof:} By $\zeta(\overline{s},y)=\overline{\zeta(s,y)}$, it 
suffices to show that 
\begin{equation}
\int\limits_{2\pi}^{V} \left\vert \zeta\left(\frac{1}{2}+it,y\right)
\right\vert^{4} {\rm d}t
\ll V(\log V)^{10}.\label{HF}
\end{equation}

By [Tch, Lemma 1], the Hurwitz zeta function satisfies an
approximate functional equation of the form
\begin{eqnarray*}
\zeta\left(\frac{1}{2}+it,y\right)
&=& \sum\limits_{1\le m\le M} (m+y)^{-1/2-it}
+\chi(1/2+it)\sum\limits_{1\le n\le N} e(-ny)n^{-1/2+it}+\\ & &
O\left(1+M^{-3/2}\vert t\vert^{1/2}\right)
\end{eqnarray*}
if $\vert t\vert\ge 2\pi$, $1\le M\le \vert t\vert$, $N\ge 1$ and  
$2\pi MN=\vert t\vert$, where $\vert \chi(1/2+it)\vert=1$. Hence, we have
\begin{eqnarray}
\int\limits_{2\pi}^{V} \left\vert \zeta\left(\frac{1}{2}+it,y\right)
\right\vert^{4} {\rm d}t 
&\ll& V+\int\limits_{2\pi}^{V} \left\vert 
\sum\limits_{1\le m\le \sqrt{t/(2\pi)}} (m+y)^{-1/2-it}
\right\vert^{4} {\rm d}t \label{H1}\\ & & 
+\int\limits_{2\pi}^{V} \left\vert 
\sum\limits_{1\le n\le \sqrt{t/(2\pi)}} e(-ny)n^{-1/2+it} 
\right\vert^{4} {\rm d}t .\nonumber
\end{eqnarray} 

By the orthogonality relation 
$$
\int\limits_{0}^{1} e(zu)\ {\rm d}u = \left\{
\begin{array}{llll} 1 & \mbox{ if } z=0,\\ \\
0 & \mbox{ if } z\in \mathbbm{Z}\setminus \{0\},\end{array} \right. 
$$
we get
\begin{equation}
\sum\limits_{1\le m\le \sqrt{t/(2\pi)}} (m+y)^{-1/2-it} =
\int\limits_{0}^{1} \sum\limits_{1\le m\le \sqrt{V}} (m+y)^{-1/2-it}
e(mu)K(t,u)\ {\rm d}u\label{H2}
\end{equation}
for $2\pi\le t\le V$, where 
$$
K(t,u):=\sum\limits_{1\le n\le \sqrt{t/(2\pi)}} e(-nu).
$$
If $2\pi\le t\le V$, then the geometric sum $K(t,u)$ can be estimated by
\begin{equation}
K(t,u)\ll \min\{\sqrt{V},\vert\vert u\vert\vert^{-1}\}.\label{H3}
\end{equation}
This yields
\begin{equation}
\int\limits_{0}^{1} \vert K(t,u)\vert\ {\rm d}u\ll \log V.\label{H4}
\end{equation}
Using H\"older$^{\prime}$s inequality, from (\ref{H2}) and
(\ref{H4}), we obtain
\begin{eqnarray}
& &
\left\vert \sum\limits_{1\le m\le \sqrt{t/(2\pi)}} (m+y)^{-1/2-it} 
\right\vert^4 \label{H5}\\ &\ll& (\log V)^3
\int\limits_{0}^{1} \left\vert 
\sum\limits_{1\le m\le \sqrt{V}} (m+y)^{-1/2-it}e(mu) \right\vert^4
\vert K(t,u)\vert\ {\rm d}u.\nonumber 
\end{eqnarray}

Employing H\"older$^{\prime}$s inequality and [Har, Lemma 3] 
after dividing the sum on the right-hand side of (\ref{H5})
into $O(\log V)$ sums of the form
$$
\sum\limits_{M<m\le 2M} (m+y)^{-1/2-it}e(mu),
$$
we obtain
\begin{equation}
\int\limits_{2\pi}^{V} \left\vert 
\sum\limits_{1\le m\le \sqrt{V}} (m+y)^{-1/2-it}e(mu)
\right\vert^{4} {\rm d}t\ll V(\log V)^6,\label{H6}
\end{equation}
where the implied $\ll$-constant does not depend on $u$. 
Combining (\ref{H3}), (\ref{H5}) and (\ref{H6}), we get
\begin{equation}
\int\limits_{2\pi}^{V} \left\vert 
\sum\limits_{1\le m\le \sqrt{t/(2\pi)}} (m+y)^{-1/2-it}
\right\vert^{4} {\rm d}t\ll V(\log V)^{10}.\label{H7}
\end{equation}
In a similar manner, we can prove 
\begin{equation}
\int\limits_{2\pi}^{V} \left\vert 
\sum\limits_{1\le n\le \sqrt{t/(2\pi)}} e(-ny)n^{-1/2+it}
\right\vert^{4} {\rm d}t\ll V(\log V)^{10}.\label{H8}
\end{equation}

Combining (\ref{H1}), (\ref{H7}) and (\ref{H8}), we obtain (\ref{HF}). 
This completes the proof.
$\Box$\\ 

To all appearances, there is no result like Theorem 4 in the literature.

We now prove Theorem 3 along the lines of the proof of [BaH, Theorem 3].
First we write
$$
F(t):=\sum\limits_{k\sim K} a_{k} (k+\theta)^{it},\ \ \ 
G(t):=\sum\limits_{l\sim L} b_{l} (l+\xi)^{it},\ \ \ 
D(t):=\sum\limits_{k\sim K} (k+\theta)^{it},
$$
$$
E(t):=\sum\limits_{l\sim L} (l+\xi)^{it}.
$$
Similar as 
in the proof of [BaH, Theorem 3], we can suppose that $K\le L\le T$, 
for otherwise 
the desired estimate follows from a classical mean value estimate for 
$G(t)$.     

Analogous to [BaH, (17)], we have
\begin{equation} 
\int\limits_0^T \vert F(\alpha t) G(\beta t) \vert^2\ {\rm d}t 
\ll (KL)^2+\log T \max\limits_{1\le V\le T} \int\limits_{V}^{2V} 
\vert D(\alpha t)E(\beta t)\vert^2\ {\rm d}t. \label{T1}
\end{equation}
We fix $V$ in the interval $1\le V\le T$ for which the maximum is 
attained.
 
In the same manner like [BaH, (19)] one can prove 
\begin{equation}
\vert D(\alpha t) \vert \ll K^{1/2} \int\limits_{-V}^V 
\left\vert \zeta\left(\frac{1}{2}+
i\sigma-i\alpha t,\theta\right)\right\vert \rho(\sigma) {\rm d}\sigma\ +\ 
\frac{K\log T}{V} \label{T2}
\end{equation}
as well as 
\begin{equation}
\vert E(\beta t) \vert \ll L^{1/2} \int\limits_{-V}^V 
\left\vert \zeta\left(\frac{1}{2}+
i\tau-i\beta t,\xi\right)\right\vert\rho(\tau) 
{\rm d}\tau\ +\ \frac{L\log T}{V}, \label{T3}
\end{equation}
where $\rho(x):=\min(1,1/\vert x\vert)$. 
Using (\ref{T2}), (\ref{T3}) and the inequality of Cauchy-Schwarz, we deduce
\begin{eqnarray}
\label{T4} & & \\ & & \int\limits_{V}^{2V} 
\vert D(\alpha t)E(\beta t)\vert^2 {\rm d}t \nonumber\\ &\ll&
\frac{(KL)^2\log^4 T}{V^3}\ +\ \frac{K^2L\log^3 T}{V^2} 
\int\limits_{-V}^V \rho(\tau) \int\limits_{V}^{2V}
\left\vert \zeta\left(\frac{1}{2}+
i\tau-i\beta t,\xi\right)\right\vert^2
 {\rm d} t\ {\rm d}\tau\ +\nonumber\\ & & 
\frac{KL^2\log^3 T}{V^2} 
 \int\limits_{-V}^V \rho(\sigma) \int\limits_{V}^{2V}
\left\vert \zeta\left(\frac{1}{2}+
i\sigma-i\alpha t,\theta\right)\right\vert^2
 {\rm d}t\ {\rm d}\sigma\ +\ KL\log^2 T \times \nonumber\\ & &
\int\limits_{-V}^V 
\int\limits_{-V}^V 
\rho(\sigma)\rho(\tau)\int\limits_{V}^{2V}\left\vert \zeta\left(\frac{1}{2}+
i\sigma-i\alpha t,\theta\right) \zeta\left(\frac{1}{2}+
i\tau-i\beta t,\xi\right)\right\vert^2 {\rm d}t\ {\rm d}\sigma 
\ {\rm d}\tau \nonumber\\ &\ll& 
\frac{(KL)^2\log^4 T}{V^3}\ +\ \frac{K^2L\log^4 T}{V^2} 
\int\limits_{-CV}^{CV} \left\vert 
\zeta\left(\frac{1}{2}+it,\xi\right)\right\vert^2
 {\rm d} t\ +\nonumber\\ & & 
\frac{KL^2\log^4 T}{V^2} 
 \int\limits_{-CV}^{CV} 
\left\vert \zeta\left(\frac{1}{2}+
it,\theta\right)\right\vert^2
 {\rm d}t\ +\ KL\log^4 T \times\nonumber\\ & &
\left(\int\limits_{-CV}^{CV} \left\vert \zeta\left(\frac{1}{2}+
it,\theta\right)\right\vert^4 {\rm d}t\right)^{1/2}  
\left(\int\limits_{-CV}^{CV} \left\vert \zeta\left(\frac{1}{2}+
it,\xi\right)\right\vert^4 {\rm d}t\right)^{1/2},\nonumber  
\end{eqnarray}
where $C$ is a certain constant which depends only on $\alpha$ and $\beta$.
From (\ref{T1}), (\ref{T4}), Theorem 4 and a similar 
{\it second} power moment estimate for the Hurwitz zeta function (which can
be derived directly from Theorem 4 using the inequality of Cauchy-Schwarz), 
we obtain 
(\ref{100}). This completes the proof of Theorem 3. $\Box$ \\ \\ \\
{\bf Acknowledgements.} \ \  This research has been supported by a Marie Curie 
Fellowship of the European Community programme ``Improving the Human Research
Potential and the Socio-Economic Knowledge Base'' under contract number
HPMF-CT-2002-02157.\\


\begin{thebibliography}{999}
\bibitem[Bai]{1} S. Baier, {\it On the $p^{\lambda}$ problem}, 
Acta Arith. 113 (2004), 77-101.
\bibitem[BaH]{2} A. Balog, G. Harman, {\it On mean values of Dirichlet
  polynomials}, Arch. Math. (Basel) 57 (1991), 581-587.
\bibitem[Har]{4} G. Harman, {\it Fractional and integral parts of
  $p^{\lambda}$}, Acta Arith. 58 (1991), 141-152.
\bibitem[Ivi]{5} A. Ivic, {\it The Riemann Zeta-Function},
Wiley-Interscience, New York 1985.
\bibitem[Tch]{6} N. Tchudakoff, {\it On Goldbach-Vinogradov's theorem}, Ann. 
Math. 48 (1947), 515-545.
\end{thebibliography}
\end{document}